# RELATIONSHIP BETWEEN NEUMANN SOLUTIONS FOR TWO-PHASE LAMÉ-CLAPEYRON-STEFAN PROBLEMS WITH CONVECTIVE AND TEMPERATURE BOUNDARY CONDITIONS

by

*Domingo Alberto TARZIA*

Departamento de Matemática, Universidad Austral, Paraguay 1950,
S2000FZF Rosario, Argentina
and CONICET, Argentina[1]

We obtain for the two-phase Lamé-Clapeyron-Stefan problem for a semi-infinite material an equivalence between the temperature and convective boundary conditions at the fixed face in the case that an inequality for the convective transfer coefficient is satisfied. Moreover, an inequality for the coefficient which characterizes the solid-liquid interface of the classical Neumann solution is also obtained. This inequality must be satisfied for data of any phase-change material, and as a consequence the result given in Tarzia, Quart. Appl. Math., 39 (1981), 491-497 is also recovered when a heat flux condition was imposed at the fixed face.

Key words: *Lamé-Clapeyron-Stefan Problem, PCM, free boundary problem, explicit solutions, similarity solutions, Neumann solution, convective boundary condition.*

## 1. Introduction

Heat transfer problems with a phase-change such as melting and freezing have been studied in the last century due to their wide scientific and technological applications [1-10]. A review of a long bibliography on moving and free boundary problems for phase-change materials (PCM) for the heat equation is shown in [11].

We consider an homogeneous semi-infinite material which has an initial constant temperature $T_i$ and at time $t = 0$ on the fixed boundary $x = 0$ we impose a constant temperature $T_0$ ($< T_f < T_i$) where $T_f$ is the phase-change temperature. These two initial and boundary conditions imply an instantaneous phase-change process. The two-phase Lamé-Clapeyron-Stefan problem [3, 12-14] is formulated in the following way: find the free boundary $x = s(t)$, defined for $t > 0$, and the temperature $T = T(x,t)$ defined by:

$$T(x,t) = \begin{cases} T_s(x,t) < T_f & if \ 0 < x < s(t), \ t > 0 \\ T_f & if \ x = s(t), \ t > 0 \\ T_\ell(x,t) > T_f & if \ s(t) < x, \ t > 0 \end{cases} \quad (1)$$

for $x > 0$ and $t > 0$, such that the following equations and conditions are satisfied (problem ($P_1$)):

$$\rho c_s T_{s_t} - k_s T_{s_{xx}} = 0, \qquad 0 < x < s(t), \quad t > 0, \quad (2)$$

---

[1] Corresponding author: e-mail: DTarzia@austral.edu.ar

$$\rho c_\ell T_{\ell_t} - k_\ell T_{\ell_{xx}} = 0, \qquad x > s(t), \quad t > 0, \tag{3}$$

$$s(0) = 0, \tag{4}$$

$$T_\ell(x,0) = T_\ell(+\infty,t) = T_i > T_f, \qquad x > 0, \quad t > 0, \tag{5}$$

$$T_s(s(t),t) = T_f, \qquad t > 0, \tag{6}$$

$$T_\ell(s(t),t) = T_f, \qquad t > 0, \tag{7}$$

$$k_s T_{s_x}(s(t),t) - k_\ell T_{\ell_x}(s(t),t) = \rho \ell \dot{s}(t), \quad t > 0, \tag{8}$$

$$T_s(0,t) = T_0 < T_f, \qquad t > 0. \tag{9}$$

The solution of the free boundary problem ($P_1$) is the classical Neumann explicit solution [3, 13, 14] given by:

$$T_s(x,t) = T_0 + \frac{T_f - T_0}{erf(\xi\sqrt{b})} erf\left(\frac{x}{2\sqrt{\alpha_s t}}\right), \quad 0 \le x \le s(t), t > 0, \tag{10}$$

$$T_\ell(x,t) = T_i - \frac{T_i - T_f}{erfc(\xi)} erfc\left(\frac{x}{2\sqrt{\alpha_\ell t}}\right), \quad s(t) \le x, \quad t > 0, \tag{11}$$

$$s(t) = 2\xi\sqrt{\alpha_\ell t}, \qquad \left(\alpha_\ell = \frac{k_\ell}{\rho c_\ell}, \quad \alpha_s = \frac{k_s}{\rho c_s}\right), \tag{12}$$

and the dimensionless coefficient $\xi > 0$ is the unique solution of the following equation:

$$G(x) = x, \quad x > 0 \tag{13}$$

with

$$G(x) = b_4 \, F_2(\sqrt{b}x) - b_3 \, F_1(x), \tag{14}$$

$$erf(x) = \frac{2}{\sqrt{\pi}} \int_0^x \exp(-u^2) du, \quad erfc(x) = 1 - erf(x), \tag{15}$$

$$F_1(x) = \frac{exp(-x^2)}{erfc(x)}, \quad F_2(x) = \frac{exp(-x^2)}{erf(x)}, \tag{16}$$

and the dimensionless parameters:

$$b = \frac{\alpha_\ell}{\alpha_s} > 0, \quad b_3 = \frac{c_\ell(T_i - T_f)}{\ell\sqrt{\pi}} > 0, \quad b_4 = \frac{k_s(T_f - T_0)}{\rho\ell\sqrt{\pi\alpha_s\alpha_\ell}} > 0. \tag{17}$$

Explicit solutions for Stefan-like problems can be found in [15-19]. A review of available analytical solutions, for a wide range of alternative boundary conditions and properties, are provided in [9].

The goal of this paper is to obtain the explicit solution of a similarity type for the solidification process for a semi-infinite PCM with a convective boundary condition at the fixed face $x = 0$ and the relationship with the Neumann solution $(10) - (12)$ corresponding to the solidification process with a temperature boundary condition at the fixed face $x = 0$. We remark that the results obtained in this paper are theoretical and they are valid for all PCMs which can be verified experimentally.

In Section II, we consider the instantaneous solidification process corresponding to the equations and conditions (2) – (8) and (18), and we show the equivalence with the solidification process (2) – (9) when the inequality (19) is satisfied for the positive coefficient $h_0$ which characterizes the transient heat transfer coefficient $h(t) = h_0 t^{-\frac{1}{2}}$ in the boundary condition (18). Moreover, an inequality for the dimensionless coefficient $\xi$ which characterizes the solid-liquid interface $s(t)$, given by (12), is also obtained. These results complement those obtained in [13, 20, 21].

## 2. Two-phase solidification process with a convective boundary condition at the fixed face x=0

In order to solve the phase-change process with a convective condition at the fixed face $x = 0$ approximate method were used, for example in [22-31]. In [32, 33] a convective condition is considered after a transformation in order to solve a free boundary problem for a nonlinear absorption model of mixed saturated-unsaturated flow with a nonlinear soil water diffusivity. In [21] the problem was analyzed and a closed-form expression for the solid-liquid interface and both temperatures were found when the heat transfer coefficient is time-dependent and proportional to $t^{-\frac{1}{2}}$. The solution is obtained graphically and it is incorrect particularly for sufficiently small heat transfer coefficient (see Corollary 2).

The goal of this paper is to give the mathematical analysis of this heat transfer problem, that is the solidification of a semi-infinite material which is initially at the liquid phase at the constant temperature $T_i$ and a convective cooling condition is imposed at the fixed boundary $x = 0$ for a time-dependent heat transfer coefficient of the type given in the condition (18).

Now, we consider the following free boundary problem with a convective boundary condition at the fixed face $x = 0$ (problem ($P_2$)): find the free boundary $x = s(t)$, defined for $t > 0$, and the temperature $T = T(x,t)$ defined as in (1) for $x > 0$ and $t > 0$, such that the equations and boundary conditions (2)-(8) and the convective boundary condition at the fixed face:

$$k_s T_{s_x}(0,t) = \frac{h_0}{\sqrt{t}}(T_s(0,t) - T_\infty), \quad t > 0 \quad (h_0 > 0),$$ (18)

are satisfied where $T_\infty$ is the bulk temperature at a large distance from the fixed face $x = 0$, $h_0$ is the coefficient which characterizes the transfer heat coefficient at the fixed face $x = 0$ and we suppose $T_\infty < T_f < T_i$ in order to guarantee a solidification process.

**Theorem 1** [20] *If the coefficient $h_0$ satisfies the inequality:*

$$h_0 > \frac{k_i}{\sqrt{\pi \alpha_\ell}} \frac{T_i - T_f}{T_f - T_\infty},$$ (19)

*then there exists an instantaneous solidification process and the free boundary problem ($P_2$) has the unique solution of a similarity type given by:*

$$T_s(x,t) = T_\infty + \frac{(T_f - T_\infty)\left(1 + \dfrac{h_0\sqrt{\pi\alpha_s}}{k_s} erf\left(\dfrac{x}{2\sqrt{\alpha_s t}}\right)\right)}{1 + \dfrac{h_0\sqrt{\pi\alpha_s}}{k_s} erf\left(\lambda\sqrt{\dfrac{\alpha_\ell}{\alpha_s}}\right)}$$

$$= T_f - \frac{h_0\sqrt{\pi\alpha_s}(T_f - T_\infty)}{k_s} \frac{erf\left(\lambda\sqrt{\dfrac{\alpha_\ell}{\alpha_s}}\right) - erf\left(\dfrac{x}{2\sqrt{\alpha_s t}}\right)}{1 + \dfrac{h_0\sqrt{\pi\alpha_s}}{k_s} erf\left(\lambda\sqrt{\dfrac{\alpha_\ell}{\alpha_s}}\right)} \qquad (20)$$

$$T_\ell(x,t) = T_i - (T_i - T_f)\frac{erfc\left(\dfrac{x}{2\sqrt{\alpha_\ell t}}\right)}{erfc(\lambda)} = T_f + (T_i - T_f)\left(1 - \frac{erfc\left(\dfrac{x}{2\sqrt{\alpha_\ell t}}\right)}{erfc(\lambda)}\right), \qquad (21)$$

$$s(t) = 2\lambda\sqrt{\alpha_\ell t}\,, \qquad (22)$$

*and the dimensionless coefficient $\lambda > 0$ satisfies the following equation:*

$$F(x) = x, \quad x > 0, \qquad (23)$$

*where the function $F$ and the parameters $b_i$ are given by:*

$$F(x) = b_1 \frac{\exp(-bx^2)}{1 + b_2 erf(x\sqrt{b})} - b_3 F_1(x), \qquad (24)$$

$$b_1 = \frac{h_0(T_f - T_\infty)}{\rho\ell\sqrt{\alpha_\ell}} > 0, \qquad b_2 = \frac{h_0}{k_s}\sqrt{\pi\alpha_s} > 0. \qquad (25)$$

**Remark 1** *If the coefficient $h_0$ satisfies the inequality (19) then for the solid and liquid temperatures, given by (20) and (21) respectively, we have the following properties:*

$$T_\infty < T_s(x,t) < T_f\,, \quad 0 < x < s(t), \quad t > 0, \qquad (26)$$

$$T_f < T_\ell(x,t) < T_i\,, \quad x > s(t), \quad t > 0. \qquad (27)$$

**Corollary 2** *If the coefficient $h_0$ satisfies the following inequalities:*

$$0 < h_0 \leq \frac{k_l}{\sqrt{\pi\alpha_\ell}} \frac{T_i - T_f}{T_f - T_\infty}. \qquad (28)$$

*then the free boundary problem $(\mathrm{P}_2)$ is a classical heat transfer problem for the initial liquid phase whose solution is given by:*

$$T_\ell(x,t) = \frac{T_\infty + \frac{k_\ell T_i}{h_0\sqrt{\alpha_\ell \pi}} + (T_i - T_\infty)erf\left(\frac{x}{2\sqrt{\alpha_\ell t}}\right)}{1 + \frac{k_\ell}{h_0\sqrt{\alpha_\ell \pi}}} = T_\infty + \frac{T_i - T_\infty}{1 + \frac{k_\ell}{h_0\sqrt{\alpha_\ell \pi}}}\left(\frac{k_\ell}{h_0\sqrt{\alpha_\ell \pi}} + erf\left(\frac{x}{2\sqrt{\alpha_\ell t}}\right)\right)$$

$$= T_i - \frac{T_i - T_\infty}{1 + \frac{k_\ell}{h_0\sqrt{\alpha_\ell \pi}}}erfc\left(\frac{x}{2\sqrt{\alpha_\ell t}}\right), \quad x > 0, \quad t > 0. \tag{29}$$

*Moreover, the temperature (29) solves the heat transfer problem without a phase-change process given by the following equation, and initial and boundary conditions:*

$$\rho c_\ell T_{\ell_t} - k_\ell T_{\ell_{xx}} = 0, \qquad x > 0, \quad t > 0, \tag{30}$$

$$T_\ell(x,0) = T_\ell(+\infty,t) = T_i > T_f, \qquad x > 0, \quad t > 0, \tag{31}$$

$$k_\ell T_{\ell_x}(0,t) = \frac{h_0}{\sqrt{t}}(T_\ell(0,t) - T_\infty), \quad t > 0. \tag{32}$$

**Proof.** The temperature, given by (29), verifies easily the heat equation (30) and the conditions (31)-(32) and it has the following properties:

$$T_\infty < T_\ell(x,t) < T_i, \quad \forall x > 0, t > 0. \tag{33}$$

Moreover, the temperature at the boundary $x = 0$ is given by:

$$T_\ell(0,t) = \frac{T_\infty + \frac{k_\ell T_i}{h_0\sqrt{\alpha_\ell \pi}}}{1 + \frac{k_\ell}{h_0\sqrt{\alpha_\ell \pi}}} = T_i - \frac{T_i - T_\infty}{1 + \frac{k_\ell}{h_0\sqrt{\alpha_\ell \pi}}} = T_\infty + \frac{\frac{k_\ell}{h_0\sqrt{\alpha_\ell \pi}}(T_i - T_\infty)}{1 + \frac{k_\ell}{h_0\sqrt{\alpha_\ell \pi}}}, \quad t > 0, \tag{34}$$

and it satisfies the following inequalities:

$$T_\infty < T_\ell(0,t) < T_i, \quad \forall t > 0. \tag{35}$$

Therefore, the temperature, given by (29), verifies the following properties:

$$T_f \le T_\ell(0,t) < T_\ell(x,t) < T_i, \quad \forall x > 0, t > 0, \tag{36}$$

owing to inequality (28) and the following inequalities:

$$T_\ell(x,t) - T_\ell(0,t) = \frac{T_i - T_\infty}{1 + \frac{k_\ell}{h_0\sqrt{\alpha_\ell \pi}}}erf\left(\frac{x}{2\sqrt{\alpha_\ell t}}\right) > 0, \quad \forall x > 0, \quad t > 0, \tag{37}$$

and

$$T_\ell(0,t) - T_f = \frac{T_i - T_\infty}{1 + \frac{h_0\sqrt{\alpha_\ell\pi}}{k_\ell}}(T_f - T_\infty) \geq \frac{T_i - T_\infty}{1 + \frac{T_i - T_f}{T_f - T_\infty}} - (T_f - T_\infty)$$

$$= (T_f - T_\infty) - (T_f - T_\infty) = 0, \quad \forall t > 0$$

(38)

□

**Remark 2** *For sufficiently small heat transfer coefficient (see inequality (28)) we have only a heat transfer problem and we can not obtain a phase-change process.*

**Theorem 3** *Let it be* $T_\infty < T_0 < T_f < T_i$. *If the coefficient* $h_0$ *satisfies the inequality (19) then the free boundary problems* ($P_1$) *and* ($P_2$) *are equivalents. Moreover, we have:*
(a) *the relationship between the datum* $T_0$ *of problem* ($P_1$) *with the data* $T_\infty$ *and* $h_0$ *of the problem* ($P_2$) *is given by:*

$$T_0 = \frac{T_f + T_\infty \frac{h_0\sqrt{\pi\alpha_s}}{k_s} erf\left(\lambda\sqrt{\frac{\alpha_\ell}{\alpha_s}}\right)}{1 + \frac{h_0\sqrt{\pi\alpha_s}}{k_s} erf\left(\lambda\sqrt{\frac{\alpha_\ell}{\alpha_s}}\right)},$$

(39)

(b) *the relationship between the data* $h_0$ *and* $T_\infty$ ($< T_0$) *of problem* ($P_2$) *with the datum* $T_0$ *of the problem* ($P_1$) *is given by:*

$$h_0 = \frac{k_s}{\sqrt{\pi\alpha_s}}\frac{T_f - T_0}{T_0 - T_\infty}\frac{1}{erf\left(\xi\sqrt{\frac{\alpha_\ell}{\alpha_s}}\right)}.$$

(40)

**Proof.** (a) If we compute the temperature at the fixed face $x = 0$ of the solution $(20) - (22)$ of the free boundary problem ($P_2$) we obtain that:

$$T_s(0,t) = T_\infty + \frac{T_f - T_\infty}{1 + \frac{h_0\sqrt{\pi\alpha_s}}{k_s} erf\left(\lambda\sqrt{\frac{\alpha_\ell}{\alpha_s}}\right)} = T_f - \frac{(T_f - T_\infty)\frac{h_0\sqrt{\pi\alpha_s}}{k_s} erf\left(\lambda\sqrt{\frac{\alpha_\ell}{\alpha_s}}\right)}{1 + \frac{h_0\sqrt{\pi\alpha_s}}{k_s} erf\left(\lambda\sqrt{\frac{\alpha_\ell}{\alpha_s}}\right)},$$

(41)

that is (33). As $T_0 < T_f$ we can consider the free boundary problem ($P_1$) whose solution is given by $(10) - (12)$ with $T_0$ defined by (39). In order to obtain the relationship between problems ($P_1$) and ($P_2$) it is sufficient to show that $\lambda = \xi$. Then, we have:

$$G(\lambda) = b_4 F_2(\sqrt{b}\lambda) - b_3 F_1(\lambda) = b_1\frac{\exp(-b\lambda^2)}{1 + b_2 erf(\lambda\sqrt{b})} - b_3 F_1(\lambda) = F(\lambda) = \lambda,$$

(42)

and by the uniqueness of the solution of the equation (13) we get $\lambda = \xi$. Then, the temperatures and solid-liquid interfaces of both problems ($P_1$) and ($P_2$) are equals.
(b) Conversely, if we consider the solution of the problem ($P_1$), we compute $T_s(0,t)$ and $T_{s_x}(0,t)$, then the positive coefficient $h_0$, which satisfies the boundary condition (18), can be defined by (40). As $h_0 > 0$ we can consider the free boundary problem ($P_2$) whose solution is given by $(20) - (22)$ with $h_0$ defined by (40) and $T_\infty < T_0$. In order to obtain the relationship between problems ($P_1$) and ($P_2$) it is sufficient to show that $\xi = \lambda$. Then, we have:

$$F(\xi) = b_1 \frac{\exp(-b\xi^2)}{1 + b_2 erf(\xi\sqrt{b})} - b_3 F_1(\xi) = b_4 F_2(\sqrt{b}\xi) - b_3 F_1(\xi) = G(\xi) = \xi , \qquad (43)$$

and by the uniqueness of the solution of the equation (23) we get $\xi = \lambda$. Then, the temperatures and solid-liquid interfaces of both problems ( $P_1$ ) and ( $P_2$ ) are equals.

Therefore, the two problems ( $P_1$ ) and ( $P_2$ ) are equivalents and the thesis holds. $\qquad\square$

**Remark 3** *Owing to the equivalence proved in Theorem 3, when we impose the two conditions (5) and (9) to a PCM we have an instantaneous phase-change process and then the corresponding heat transfer coefficent at the fixed face* $x = 0$ *is time-dependent and proportional to* $t^{-\frac{1}{2}}$ *of the type (18) and the coefficient* $h_0$ *must verify the inequality (19).*

**Theorem 4** *Let it be* $T_\infty < T_0 < T_f < T_i$ . *The coefficient* $\xi = \xi_{T_0}$ *which characterizes the free boundary (12) of the Neumann solution of the problem ( $P_1$ ) (with datum* $T_0 (<T_f)$ *at the fixed face* $x = 0$ *) satisfies the following inequality:*

$$erf\left(\xi\sqrt{\frac{\alpha_\ell}{\alpha_s}}\right) < \frac{k_s}{k_\ell}\sqrt{\frac{\alpha_\ell}{\alpha_s}}\frac{T_i - T_\infty}{T_0 - T_\infty}\frac{T_f - T_0}{T_i - T_f}, \quad \forall T_\infty < T_0 . \qquad (44)$$

**Proof.** By Theorem 3, we have that $h_0$, defined by (40), must satisfy the inequality (19) and therefore the coefficient $\xi$ of the solid-liquid interface of the Neumann solution of the problem ( $P_1$ ) must also satisfy the inequality (44). $\qquad\square$

**Remark 4** *The inequality (44) has a physical meaning for the classical Neumann solution (10) – (12) when the parameters of the problem ( $P_2$ ) verify the inequality:*

$$\frac{k_s}{k_\ell}\sqrt{\frac{\alpha_\ell}{\alpha_s}}\frac{T_i - T_\infty}{T_0 - T_\infty}\frac{T_f - T_0}{T_i - T_f} < 1, \quad \text{with } T_\infty < T_0 < T_f < T_i . \qquad (45)$$

**Corollary 5** *The coefficient* $\xi = \xi_{T_0}$ *which characterizes the free boundary (12) of the Neumann solution of the problem ( $P_1$ ) (with* $T_0 (<T_f < T_i)$ *at the fixed face* $x = 0$ *) also satisfies the following inequality:*

$$erf\left(\xi\sqrt{\frac{\alpha_\ell}{\alpha_s}}\right) < \frac{k_s}{k_\ell}\sqrt{\frac{\alpha_\ell}{\alpha_s}}\frac{T_f - T_0}{T_i - T_f} . \qquad (46)$$

**Proof.** The function of the right hand side in (44) is a strictly increasing function of the variable $T_\infty$. Then, the inequality (46) is obtained by taking in the inequality (38) the limit when $T_\infty \to -\infty$. $\qquad\square$

**Remark 5** *The inequality (46) is the same inequaliy as obtained in* [13, 34] *where a two-phase Lamé-Clapeyron-Stefan free boundary problem with a heat flux of the type*

$$k_s T_{s_x}(0,t) = q(t), t > 0 \quad \left(q(t) = \frac{q_0}{\sqrt{t}}, q_0 > 0\right), \qquad (47)$$

*at the fixed face* $x = 0$ *is imposed instead of the temperature boundary condition (9).*

**Theorem 6** *Let it be* $T_\infty < T_0 < T_f < T_i$ . *The coefficient* $\lambda = \lambda(h_0)$ *of the free boundary (22) is a strictly increasing function of the variable* $h_0$, *defined in the interval* $\left( k_\ell (T_i - T_f)(T_f - T_\infty)^{-1}(\pi \alpha_\ell)^{-1/2}, +\infty \right)$, *and it satisfies the following properties:*

$$0 = \lambda \left( \frac{k_\ell}{\sqrt{\pi \alpha_\ell}} \frac{T_i - T_f}{T_f - T_\infty} \right) < \lambda = \lambda(h_0) < \lambda_{T_\infty}, \quad \forall h_0 > \frac{k_\ell}{\sqrt{\pi \alpha_\ell}} \frac{T_i - T_f}{T_f - T_\infty}, \qquad (48)$$

*where* $\lambda_{T_\infty}$ *is the coefficient which characterizes the Neumann solution of the free boundary problem (2)-(8) and the following temperature boundary condition at the fixed face* $x = 0$ *given by:*

$$T_s(0,t) = T_\infty, t > 0, \qquad (49)$$

*instead of the boundary condition (9).*

**Proof.** The positive coefficient $\lambda$ is defined as the unique solution of the equation (23). If we consider $\lambda = \lambda(h_0)$ as a function of the parameter $h_0$ which must verify the inequality (19), we can obtain that $b_1 = b_1(h_0)$ and $b_2 = b_2(h_0)$, and then $F = F(h_0)$ are strictly increasing functions of the parameter $h_0$, therefore $\lambda = \lambda(h_0)$ is also a strictly increasing function of the variable $h_0$.

If $h_0 = \frac{k_\ell}{\sqrt{\pi \alpha_\ell}} \frac{T_i - T_f}{T_f - T_\infty}$ then $b_1 = b_3$ and so we have $\lambda \left( \frac{k_\ell}{\sqrt{\pi \alpha_\ell}} \frac{T_i - T_f}{T_f - T_\infty} \right) = 0$. By the other

part, when we consider the limit $h_0 \to \infty$ in equation (23) we obtain that

$$\lambda(+\infty) = \lambda_{T_\infty} \qquad (50)$$

where $\lambda_{T_\infty}$ is the solution of the equation (13) when we impose the boundary condition (49) at the fixed face $x = 0$ instead of the boundary condition (9). $\quad \square$

**Remark 6** *In [35] an equivalence of two problems for the one-phase fractional Lamé-Clapeyron-Stefan problems has been considered. For the classical one-phase problems the equivalence is trivial because* $T_i = T_f$ *and therefore the inequality (19) is always verified*

$$h_0 > \frac{k_\ell}{\sqrt{\pi \alpha_\ell}} \frac{T_i - T_f}{T_f - T_\infty} = 0 \quad . \qquad (51)$$

*Following [35, 36] we can also generalize the inequality (19) for the two-phase fractional Lamé-Clapeyron-Stefan problem in a forthcoming paper [37].*

**Remark 7** *At the last, for a suggestion of an anonymous referee, we will transform the problem* $(P_2)$, *given by the equations and conditions (2)-(8) and (18), and the inequality (19) in a dimensionless form. We define the following dimensionless change of variables:*

$$\eta = \frac{x}{L}, \quad \tau = \frac{\alpha_s t}{L^2}, \quad r(\tau) = \frac{s(t)}{L}, \quad \theta_\ell(\eta, \tau) = \frac{T_\ell(x,t) - T_f}{T_i - T_f}, \quad \theta_s(\eta, \tau) = \frac{T_s(x,t) - T_f}{T_i - T_f} \qquad (52)$$

*where* $L$ *is a characteristic length. Therefore, the equations and conditions (2)-(8) and (18) are transformed as:*

$$\theta_{s_\tau} - \theta_{s_{\eta\eta}} = 0, \qquad 0 < \eta < r(\tau), \quad \tau > 0, \tag{53}$$

$$\theta_{\ell_\tau} - \frac{\alpha_\ell}{\alpha_s} \theta_{\ell_{\eta\eta}} = 0, \qquad \eta > r(\tau), \quad \tau > 0, \tag{54}$$

$$r(0) = 0, \tag{55}$$

$$\theta_\ell(\eta, 0) = \theta_\ell(+\infty, \tau) = 1, \qquad \eta > 0, \quad \tau > 0, \tag{56}$$

$$\theta_s(r(\tau), \tau) = 0, \qquad \tau > 0, \tag{57}$$

$$\theta_\ell(r(\tau), \tau) = 0, \qquad \tau > 0, \tag{58}$$

$$\theta_{s_\eta}(r(\tau), \tau) - \frac{k_\ell}{k_s} \theta_{\ell_\eta}(r(\tau), \tau) = \frac{1}{Ste} r'(\tau), \quad \tau > 0, \tag{59}$$

$$\theta_{s_\eta}(0, \tau) = \frac{B}{\sqrt{\tau}} \left( \theta_s(0, \tau) + \theta_\infty \right), \quad \tau > 0, \tag{60}$$

where $Ste$ is the Stefan number and $B/\sqrt{\tau}$ is the Biot number defined by the following expressions:

$$Ste = \frac{c_s(T_i - T_f)}{\ell} > 0, \tag{61}$$

$$B = \frac{h_0 \sqrt{\alpha_s}}{k_s} > 0, \tag{62}$$

and

$$\theta_\infty = \frac{T_f - T_\infty}{T_i - T_f} > 0. \tag{63}$$

Moreover, the inequality (19) for the coefficient $h_0$, which characterized the heat transfer coefficient in the boundary condition (18), is transformed in the following way:

$$B > \frac{k_l}{k_s} \sqrt{\frac{\alpha_s}{\pi \alpha_\ell}} \frac{1}{\theta_\infty} = \frac{k_l}{k_s} \sqrt{\frac{\alpha_s}{\pi \alpha_\ell}} \frac{T_i - T_f}{T_f - T_\infty}. \tag{64}$$

In the case we define the dimensionless time as $\tau = \frac{\alpha_l t}{L^2}$ we can obtain similar results as before but with small differences in the equations (53) and (54), and conditions and definitions (59), (61), (62) and (64) respectively.

## Conclusions

An equivalence between two Lamé-Clapeyron-Stefan problems with a temperature and a convective boundary conditions at the fixed face of a semi-infinite PCM is obtained for sufficiently large heat transfer coefficient. Then, an inequality for the coefficient which characterizes the solid-liquid interface of the classical Neumann solution is also obtained.

**Aknowledgements**

We would like to thank three anonymous referees for their constructive comments which improved the readability of the manuscript. The present work has been partially sponsored by the Project PIP No 0534 from CONICET-UA, Rosario, Argentina, and Grant AFORS FA9550-14-1-0122.

**Nomenclature**

$b, b_i$    Parameters defined in (17) and (25) (i=1, 2, 3, 4), [-]

$c$    Specific heat, $[\,J\,kg^{-1}K^{-1}\,]$

$h$    Transient heat transfer coefficient at x=0, $[\,kg\,K^{-1}s^{-3}\,]$

$h_o$    Coefficient that characterizes the heat transfer coefficient at x=0, $[\,kg\,K^{-1}s^{-5/2}\,]$

$k$    Thermal conductivity, $[\,Wm^{-1}K^{-1}\,]$

$\ell$    Latent heat of fusion by unit of mass, $[\,J\,kg^{-1}\,]$

$q$    Transient heat flux at x=0, $[\,kg\,s^{-3}\,]$

$q_o$    Coefficient that characterizes the transient heat flux at x=0, $[\,kg\,s^{-5/2}\,]$

$s$    Position of the solid-liquid interface, [m]

$t$    Time, [s]

$T$    Temperature, [ºC]

$T_f$    Phase-change temperature $(T_0 < T_f < T_i)$, [ºC]

$T_\infty$    Bulk temperature at the fixed face $x = 0$ $(T_\infty < T_f < T_i)$, [ºC]

$x$    Spatial coordinate, [m]

Greek symbols

$\alpha$    Diffusivity coefficient $\left(= k(\rho c)^{-1}\right)$, $[\,m^2 s^{-1}\,]$

$\lambda$    Coefficient that characterizes the free boundary in Eq. (22), [-]

$\rho$    Density, $[\,kg\,m^{-3}\,]$

$\xi$    Coefficient that characterizes the free boundary in Eq. (12), [-]

Subscripts

$l$    Liquid phase

$s$    Solid phase